# Augmented Zagreb Index of Polyhex Nanotubes

Nazeran Idrees[1], Afshan Sadiq[2], Muhammad Jawwad Saif[3]

[1]Department of Mathematics, Government College University Faisalabad, Pakista
[2]Abdus Salam School of Mathematical Sciences, Government College University Lahore, Pakistan
[3]Department of Applied Chemistry, Government College University Faisalabad, Pakistan

*email: corresponding author*: nazeranjawwad@gmail.com (Nazeran Idrees),           2. afshansadiq6@gmail.com, 3. jawwadsaif@gmail.com,

## ABSTRACT

Augmented Zagreb Index is a newly defined degree based topological invariant which has been well established for its better correlation properties and is defined as $AZI(G) = \sum_{uv \in E(G)} \left( \frac{d_G(u)d_G(v)}{d_G(u) + d_G(v) - 2} \right)^3$, where $E(G)$ is the edge set of graph $G$ and $d(u), d(v)$ are the degrees of the end vertices $u$ and $v$ of edge $uv$, respectively. It has outperformed many well known degree based topological indices. In this article we give closed formulae for the augmented Zagreb index of arm-chair polyhex and zigzag edge polyhex nanotubes.

## 1. Introduction

A molecular graph is a representation of a chemical compound having atoms as vertices and the bonds between atoms correspond to the edge of the graph. The collection of vertices of a graph, say $G$, is denoted by $V(G)$ and the set of edges is denoted by $E(G)$. The degree of a vertex $v$ of a graph is the number of vertices of $G$ adjacent to $v$, denoted by $d_G(v)$ or simply as $d_v$.

Topological invariant of a graph is a single number descriptor which is correlated to certain chemical, thermodynamical and biological behaviour of the chemical compounds. Several topological indices have been defined over past decades which depend on degree of vertices and the Randic index is the very first of these indices, introduced by Milan Randic,[1] formulated as

$$R(G) = \sum_{uv \in E(G)} \frac{1}{\sqrt{d_G(u)d_G(v)}} \qquad (1)$$

Atom bond connectivity index (ABC), introduced by Estrada et al.,[2] defined as

$$ABC(G) = \sum_{uv \in E(G)} \sqrt{\frac{d_G(u) + d_G(v) - 2}{d_G(u)d_G(v)}} \qquad (2)$$

Recently Furtula et al,[3,4] introduced "augmented Zagrab index", denoted by AZI and is defined as

$$AZI(G) = \sum_{uv \in E(G)} \left( \frac{d_G(u)d_G(v)}{d_G(u) + d_G(v) - 2} \right)^3 \qquad (3)$$

Some basic investigation implied that AZI index has better correlation properties and structural sensitivity among the very well established degree based topological indices.[5,6]

Several topological indices have been studied for Polyhex nanotubes.[7-9] We give an explicit relation to find the AZI index of armchair polyhex and zigzag edge polyhex nanotubes.

## 2.    MAIN RESULTS AND DISCUSSION

### 2.1. Armchair Polyhex Nanotubes.
Consider the armchair polyhex nanotubes $G = TUAC_6[m,n]$, where $m$ denotes number of hexagons in first row and $n$ denotes the number of rows. The number of vertices/atoms of armchair polyhex nanotube is equal to $|V(TUAC_6[m,n])| = 2m(n+2)$ and the number of edges/bonds is $|E(TUAC_6[m,n])| = 3mn + 4m$.

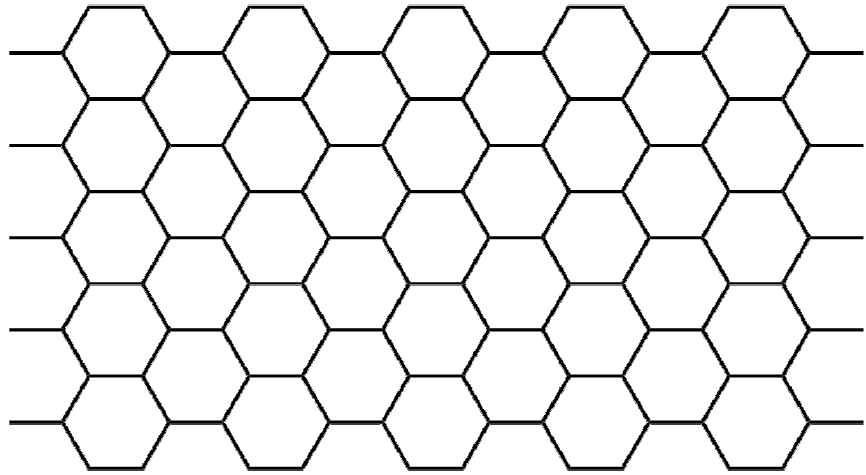

FIGURE 1. Graph of arm-chair polyhex $TUAC_6[5,9]$ nanotube.

There are three different kind of edges of $G$ depending on the degree of terminal vertices of edges, as shown in table.

TABLE 1. edge partition of 2-dimentional graph of $TUAC_6$ nanotube with respect to degree of end vertices of edges.

| Number of edges | $(d_u d_v)$ where $uv \in E(G)$ is edge |
|---|---|
| $2m$ | $(2,2)$ |
| $4m$ | $(2,3)$ |
| $3mn - 2m$ | $(3,3)$ |

**Theorem 2.1.** *Consider the $TUAC_6[m,n]$ nanotube, then its augmented Zagreb index (AZI) is equal to*

$$AZI(TUAC_6[m,n]) = \frac{2187}{64}mn - \frac{573}{64}m \qquad (4)$$

**Proof.** Consider the graph $G$ of $TUAC_6[m,n]$nanotube. The number of vertices in $TUAC_6[m,n]$ are $2m(n+2)$ and the number of edges of the nanotube is $3mn + 4m$. Now using different type of edges corresponding to the degrees of terminal vertices of edges of $G$ given in Table 1 we compute the augmented Zagreb index of $G$ which is expressed as

$$AZI(G) = \sum_{uv \in E(G)} \left( \frac{d_G(u)d_G(v)}{d_G(u) + d_G(v) - 2} \right)^3$$

This implies that

$$AZI(TUAC_6[m,n]) = 2m \left( \frac{2.2}{2 + 2 - 2} \right)^3 + 4m \left( \frac{2.3}{2 + 3 - 2} \right)^3 + (3mn - 2m) \left( \frac{3.3}{3 + 3 - 2} \right)^3$$

$$= 48m + (3mn - 2m) \left( \frac{9}{4} \right)^3$$

Which can be further simplified to

$$AZI(TUAC_6[m,n]) = \frac{2187}{64}mn - \frac{807}{32}m \approx 34mn - 25m.$$

**2.2. Zigzag-edge polyhex nanotubes.** Let $H$ be the two dimensional graph of $TUZC_6[m,n]$ nanotube, where $m$ denotes number of hexagons in first row and $n$ denotes number of repetitions.

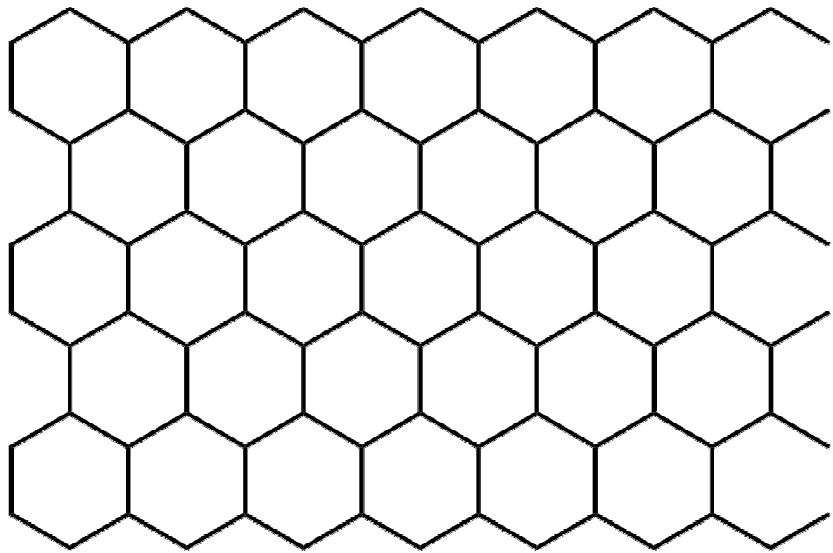

Figure 2. Graph of zigzag edge polyhex $TUZC_6[7,5]$ nanotube.

The total number of vertices $V$ of $H$ is $2mn + 2m$ and the number of edges $E$ of $H$ is equal to $3mn + 2m$. There are two different kind of edges of $H$ depending on the degree of terminal vertices of edges.

TABLE 2. edge partition of 2-dimensional graph of $TUZC_6[m,n]$ nanotube with respect to degree of end vertices of all the edges.

| Number of edges | $(d_u d_v)$ where $uv \in E(H)$ is edge |
|:---:|:---:|
| $4m$ | (2,3) |
| $3mn - 2m$ | (3,3) |

**Theorem 2.2**. *Consider the $TUZC_6[m,n]$ nanotube, then its augmented Zagreb index(AZI) is*

$$AZI(TUZC_6[m,n]) = \frac{2187}{64}mn - \frac{597}{64}m. \qquad (5)$$

Proof. Let $H$ be the graph of $TUZC_6[m,n]$ nanotube having number of vertices equal to $2mn + 2m$ and the number of edges is $3mn + 2m$. Now we compute the augmented Zagreb index of $H$ by employing the edge partition corresponding to the degree of end vertices of edges of $H$ given in Table 2. Since

$$AZI(H) = \sum_{uv \in E(G)} \left( \frac{d_G(u) d_G(v)}{d_G(u) + d_G(v) - 2} \right)^3$$

This implies that

$$AZI(TUZC_6[m,n]) = 4m\left(\frac{2.3}{2+3-2}\right)^3 + (3mn - 2m)\left(\frac{3.3}{3+3-2}\right)^3$$

After some simplification steps, we get

$$AZI(TUZC_6[m,n]) = \frac{2187}{64}mn - \frac{434}{64}m$$

$$\approx 34mn - 7m.$$

## 3. CONCLUSION

In the present study we computed the augmented Zagreb index of armchair as well as zigzag edge nanotubes and gave closed formulae of augmented Zagreb index for these nanotubes which have not been computed earlier according to the best of authors' knowledge. It is quite useful to investigate more about nanomaterials that have become the core interest of basic sciences. The AZI index can help much better to correlate the graphs with the physico chemical properties of these materials.